\documentclass{article}

\usepackage{amsfonts,amsmath,amsthm,amssymb,mathrsfs}
\usepackage{fullpage}
\usepackage{enumitem}
\usepackage[latin1]{inputenc}
\usepackage{graphicx,xcolor}
\usepackage{hyperref}
\usepackage{fancyhdr}
\usepackage[font={small,sf}]{caption}

\newtheorem{theorem}{Theorem}

\title{Enumeration of rooted 3-connected bipartite planar maps}

\author{
	Marc Noy
	\thanks{
		Departament de Matem\`atiques and Institut de Matem\`atiques (IMTech) de la Universitat Polit\`ecnica de Catalunya (UPC), and Centre de Recerca Matem\`atica (CRM), Barcelona, Spain.
		E-mail: {\tt marc.noy@upc.edu}.
  }
	\and
	Cl\'ement Requil\'e	
	\thanks{
		Departament de Matem\`atiques and Institut de Matem\`atiques (IMTech) de la Universitat Polit\`ecnica de Catalunya (UPC), Barcelona, Spain.
  	E-mail: {\tt clement.requile@upc.edu}.
	}
	\and
	Juanjo Ru\'e
	\thanks{
		Departament de Matem\`atiques and Institut de Matem\`atiques (IMTech) de la Universitat Polit\`ecnica de Catalunya (UPC), and Centre de Recerca Matem\`atica (CRM), Barcelona, Spain.
		E-mail: {\tt juan.jose.rue@upc.edu}.
  }
}

\begin{document}

\maketitle

\paragraph{Abstract:}
We provide the first solution to the problem of counting rooted 3-connected bipartite planar maps.
Our  starting point is the enumeration of  bicoloured  planar maps according to the number of edges and  monochromatic edges, following Bernardi and Bousquet-M\'elou [J. Comb. Theory Ser. B,  101 (2011), 315--377].
The decomposition of a map into 2- and 3-connected components allows us to obtain the generating functions of 2- and  3-connected bicoloured maps.
Setting to zero the variable marking monochromatic edges we obtain the  generating function of  3-connected bipartite maps, which is algebraic of degree 26.
We deduce from it an asymptotic estimate for the number of 3-connected bipartite planar maps of the form $t\, n^{-5/2}\gamma^n$, where $\gamma = \rho^{-1} \approx 2.40958$ and $\rho \approx 0.41501$ is an algebraic number of degree $10$.

\paragraph{R\'esum\'e:}
Nous apportons la premi\`ere solution au probl\`eme du d\'enombrement des cartes enracin\'ees plan-aires qui sont biparties et 3-connexes.
Notre point de d\'epart est l'\'enum\'eration des cartes planaires bi-colori\'ees, d'apr\`es Bernardi et Bousquet-M\'elou [J. Comb. Theory Ser. B,  101 (2011), 315--377].
La d\'ecomposition d'une carte en composantes 2- et 3-connexes nous permet ensuite d'obtenir les fonctions g\'en\'eratrices des cartes bi-colori\'ees 2- et 3-connexes.
En \'evaluant \`a z\'ero la variable marquant le nombre d'ar\^etes monochromes, nous obtenons alors la fonction g\'en\'eratrice des cartes biparties 3-connexes.
Cette derni\`ere est alg\'ebrique de degr\'e 26.
Nous en d\'eduisons une estimation asymptotique de la forme $t\, n^{-5/2}\gamma^n$ du nombre de cartes planaires biparties 3-connexes, avec $\gamma = \rho^{-1} \approx 2.40958$ et o\`u $\rho \approx 0.41501$ est un nombre alg\'ebrique de degr\'e $10$.

%%%%%%%%%%%%%%%%%%%%%%%%%%%%%%%%%%%%%%%%%%%%%%%%%%%%%%%%%%%%%%%%%%%%%%%%%%%%%%%%
%
%	Introduction
%
%%%%%%%%%%%%%%%%%%%%%%%%%%%%%%%%%%%%%%%%%%%%%%%%%%%%%%%%%%%%%%%%%%%%%%%%%%%%%%%%
\section{Introduction}\label{sec:introduction}

The theory of map enumeration was initiated by William Tutte in the 1960's, motivated by the Four Colour Problem, the most notorious open problem in graph theory at that time.
In his own words~\cite[Chapter 10]{T98}:

\begin{quote}
	From time to time in a graph-theoretical career one's thoughts turn to the Four Colour Problem.
	It occurred to me once that it might be possible to get results of interest in the theory of map-colourings without actually solving the Problem.
	For example it might be possible to find the average number of 4-colourings, on vertices, for planar triangulations of a given size.
	One would determine the number of triangulations of $2n$ faces, and then the number of 4-coloured triangulations of $2n$ faces.
	Then one would divide the second number by the first to get the required average.
	I gathered that this sort of retreat from a difficult problem to a related average was not unknown in other branches of Mathematics, and that it was particularly common in Number Theory.
\end{quote}

\noindent
The first task for Tutte was to decide what sort of triangulations to study.
From the point of view of colouring the natural choice were 3-connected triangulations, which correspond to maximal planar graphs.
He started drawing them and, again quoting~\cite{T98}:
\begin{quote}
	Having made no progress with the enumeration of these diagrams I bethought myself of Cayley's work on the enumeration of trees.
	His first successes had been with the rooted trees, in which one vertex is distinguished as the ``root''.
	Perhaps I should root triangulations in some way and try to enumerate the rooted ones.
\end{quote}

\noindent
Tutte defined rooted maps and succeeded in counting rooted triangulations, thus starting his famous series of `census' papers on map enumeration~\cite{T62}.
Later he counted all rooted maps as well as the 2-connected and 3-connected ones~\cite{T63}, and also Eulerian maps which by duality are in bijection with bipartite maps.
From here it is not difficult to count those which are 2-connected and bipartite~\cite{LW04}.
But counting 3-connected bipartite maps has remained an open problem.

After these preliminaries let us properly define our objects.
A planar map is a connected multigraph with a given embedding in the sphere.
All maps in this paper are rooted at a directed edge.
Maps are counted according to the number of edges and up to orientation-preserving homeomorphisms of the sphere; see  for instance~\cite{S15} for basic concepts on planar maps.
A map is said to be \textit{2-connected} if it has at least two edges, no loops, and no cut vertices; the smallest one is the digon, the map with two vertices linked by a double edge.
It is furthermore said to be \textit{3-connected} if it has at least six edges, no double edges, and no vertex separators of size two; the smallest one is $K_4$, the map of the complete graph on four vertices.
Tutte showed that the number of (rooted) maps with $n$ edges is given by
\begin{equation}\label{eq:Mn}
	A_n = \frac{2\cdot 3^n}{(n + 1)(n + 2)}\binom{2n}{n}.
\end{equation}
Notice that $A_0 = 1$, corresponding to the empty map, i.e. with one vertex and no edge.
This surprisingly simple formula was explained much later by Schaeffer~\cite{S98} in his Ph.D. thesis through a remarkable bijection with `decorated' trees, drawing on previous work by Cori and Vauquelin~\cite{CV81}.
This opened the way to the study of the asymptotic metric properties of maps and their connection with Brownian motion~\cite{CS04}, culminating in the construction of the so-called Brownian map (see~\cite{LG19} for an overview).

Tutte originally proved Formula~\eqref{eq:Mn} using a correspondance with bicubic maps (defined later).
But it can be proved more directly as follows.
For a map $\mathsf{m}$ and an edge of $\mathsf{m}$ we denote by $\mathsf{m}-e$ the supression of $e$ in $\mathsf{m}$. Then, $\mathsf{m}$ is either the empty map (namely, with one vertex and no edges), or $\mathsf{m}-e$ is the union of two disjoint maps (which can be rooted canonically and determine $\mathsf{m}$ uniquely), or $\mathsf{m}-e$ is connected.
In the latter case $\mathsf{m}-e$ can also be canonically rooted, but to recover $\mathsf{m}$ we need to know  which vertex of the root face is the second vertex of the new root edge.
Hence, one needs to refine the counting by considering the number $A_{n,k}$ of rooted maps with $n$ edges and root face of degree $k$.
This leads to a quadratic equation satisfied by the counting series (or generating function) $A(z,y)= \sum_{n,k} A_{n,k} z^n y^k$
\begin{equation}\label{eq:induction_maps}
	A(z,y) = 1 + y^2zA(z,y)^2 + yz\frac{yA(z,y) - A(z,1)}{y - 1}.
\end{equation}
In this context $y$ is called a \textit{catalytic} variable.
Notice that we cannot solve~\eqref{eq:induction_maps} directly, computing $A(z) = \sum_{n}A_n z^n = A(z,1)$ as a specialisation of $A(z,y)$.
Tutte already encountered this problem in his pioneering work on the enumeration of triangulations~\cite{T62}, and it was solved completely in a systematic way by Brown in~\cite{Brown65}.
Brown's paper introduced the so-called \emph{quadratic method} (see~\cite{BMJ06} for a far-reaching generalisation to arbitrary polynomial equations with divided differences).
Applying this method one obtains the explicit formula
\begin{equation*}
	A(z) = \frac{18z - 1 + (1 - 12z)^{3/2}}{54z^2}.
\end{equation*}
In particular, $A(z)$ is an algebraic function of degree two.
From the former expression, it is straightforward to obtain~\eqref{eq:Mn} and the estimate
\begin{equation*}
	[z^n]A(z) \sim \frac{2}{\sqrt{\pi}} \, n^{-5/2} \, 12^n
	\qquad\text{as } n\to\infty,
\end{equation*}
where we use the notation $[z^n]A(z)$ to denote the $n$-th coefficient of the power series $A(z)$.
The exponent $-5/2$ is universal for `naturally defined' classes of planar maps and has been explained in a number of different ways (see for instance~\cite{DNY22} in relation to the quadratic method).

\begin{figure}
	\begin{center}
			\raisebox{.1cm}{\includegraphics[scale=1]{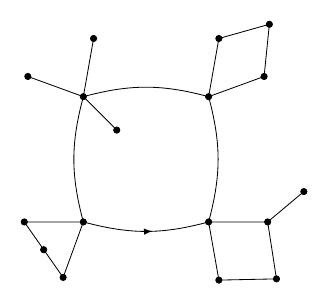}}
			\qquad\qquad\qquad
			\includegraphics[scale=1]{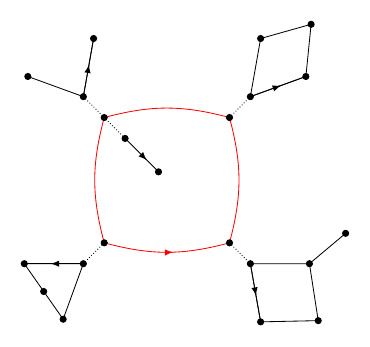}
	\end{center}
	\caption{
	The decomposition of a rooted planar map $\textsf{m}$ into its 2-core $\mathsf{C}(\textsf{m})$ (in red), i.e. the non-separable map containing the root, and maps attached by their root vertex to the corners of $\mathsf{C}(\textsf{m})$.
	}\label{fig:non-separable_map}
\end{figure}

Once $A(z)$ is determined, one can obtain the generating functions  $B(z)$ and $T(z)$ counting 2-connected and 3-connected planar maps, respectively, where $z$ marks edges.
A \textit{block} in a rooted map $\mathsf{m}$ is a maximal component, in the sense of submap,  that is either the single edge or is 2-connected.
If the unique maximal block of $\mathsf{m}$ containing the root edge is 2-connected, it is called the \textit{2-core} of $\mathsf{m}$ and is denoted by $\mathsf{C}(\mathsf{m})$.
Then $\mathsf{m}$ can be recovered by placing a rooted map at each \textit{corner} of $\mathsf{C}(\mathsf{m})$, i.e. identifying the root vertex of the rooted map with the vertex of the corner, where a corner consists of two consecutive edges with a common vertex.
See Figure~\ref{fig:non-separable_map} for an example of such decomposition.
Taking into account that the number of corners in a map is twice the number of edges, one obtains the relation
\begin{equation}\label{eq:MN}
	A(z) = 1 + 2zA(z)^2 + B(zA(z)^2),
\end{equation}
where the term 1 encodes the empty map, $2zA(z)^2$ encodes the corner substitution on the single loop and edge maps, and $B(zA(z)^2)$ encodes the corner substitution into the 2-core.
By inverting this equation one obtains
\begin{equation*}
	B(z) = z^2 + 2z^3 + 6z^4 + 22z^5 + \dots,
\end{equation*}
where the term $z^2$ corresponds to the digon.

\begin{figure}[h]
	\begin{center}
		\begin{minipage}{.45\textwidth}
		\centering
			\includegraphics[scale=1]{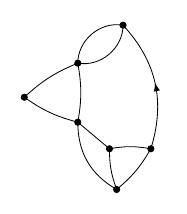}
			\quad
			\includegraphics[scale=1]{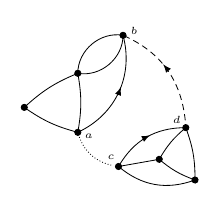}
		\end{minipage}
		\hfill\vline\hfill
		\begin{minipage}{.48\textwidth}
		\centering
			\raisebox{.2cm}{\includegraphics[scale=1]{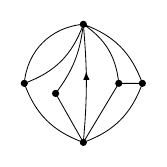}}
			\quad
			\includegraphics[scale=1]{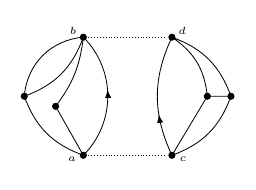}
		\end{minipage}
	\end{center}
	\caption{
	Decomposition of a series map (left) and its dual parallel map (right) into smaller 2-connected maps.
	The series map is obtained from the two smaller maps by removing the two root edges, identifying vertex $a$ with $c$, and adding a root edge from $d$ to $b$.
	The (dual) parallel map is obtained from the two smaller maps by identifying vertex $a$ with $c$ and vertex $b$ with $d$ then removing the edge $cd$.
	}\label{fig:series_parallel_maps}
\end{figure}

\begin{figure}
	\begin{center}
			\includegraphics[scale=1]{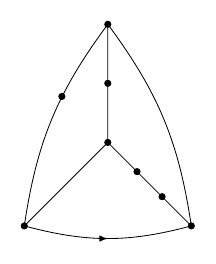}
			\qquad\qquad\qquad
			\includegraphics[scale=1]{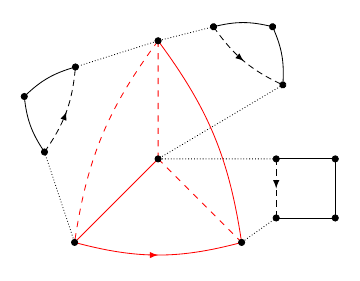}
	\end{center}
	\caption{
	Decomposition of a polyhedral map $\textsf{m}$ into its 3-core $\mathsf{T}(\textsf{m})$ (in red) and smaller 2-connected maps.
	}\label{fig:polyhedral_map}
\end{figure}

In order to find $T(z)$, one uses the decomposition of 2-connected graphs into  3-connected components (see~\cite{T58} and~\cite{T63} for details).
Let $D(z) = B(z)/z$, meaning that the root edge is not counted.
Then
\begin{equation}\label{eq:BT}
	D(z) = z + S(z) + P(z) + \frac{T(D(z))}{D(z)},
\end{equation}
where $S(z)$ and $P(z)$ encode \textit{series} and \textit{parallel} maps, respectively.
A series map is a 2-connected map obtained by the series composition of an arbitrary 2-connected map, or the single edge map, and a non-series 2-connected map, or the single edge map.
We then have
\begin{equation*}
	S(z) = D(z)(D(z) - S(z)),  \qquad P(z) =S(z), %D(z)(D(z) - P(z)).
\end{equation*}
where the second equation follows by duality.
In Figure~\ref{fig:series_parallel_maps} are examples of series and parallel decompositions of maps.
The 2-connected maps $\textsf{m}$ counted by the term $T(D(z))/D(z)$ in~\eqref{eq:BT} are called \textit{polyhedral} and are obtained after substituting 2-connected maps for the non-root edges of a 3-connected map $\textsf{T}(\textsf{m})$.
An edge $e$ of $\mathsf{T}(\textsf{m})$ is replaced by a 2-connected map $\textsf{m}'$ with root edge $r$ by identifying the two vertices of $r$ with the two vertices of $e$, then deleting both $e$ and $r$ from the resulting map.
In that case, $\textsf{T}(\textsf{m})$ is the maximal 3-connected component of $\textsf{m}$ containing the root and is called the \textit{3-core} of $\textsf{m}$.
An illustration can be found in Figure~\ref{fig:polyhedral_map}.
It follows that
\begin{equation*}
	\frac{T(z)}{z} = z - \frac{2z^2}{1 + z} - F(z),
\end{equation*}
where $F(z)$ is the functional inverse of $D(z)$.
From this expression, one obtains a quadratic equation satisfied by $T(z)$, which can be solved directly.
Using the previous algebraic expressions and by means of algebraic inversion, Tutte obtained implicit expressions for $B(z)$ and $T(z)$, which are algebraic series of degree three and two, respectively,  and deduced the estimates
\begin{equation*}
	[z^n]B(z) \sim \frac{2\sqrt{3}}{27\sqrt{\pi}} \, n^{-5/2} \left(\frac{27}{4}\right)^n
	\quad\text{and}\quad
	[z^n]T(z) \sim \frac{2}{243\sqrt{\pi}} \, n^{-5/2} \, 4^n,
	\qquad\text{as } n\to\infty.
\end{equation*}

\paragraph{From uncolored to bipartite maps.}

Our aim in this paper is to carry out the same program for \textit{bipartite} maps.
By convention we will always assume that   vertices are coloured black and white, and that the root vertex is always  black.
A straightforward modification of the quadratic method (see for instance Equation (3.3) in~\cite{BBM11}) provides the series of non-empty bipartite maps as
\begin{equation}\label{GF:Mb}
		A_b(z) = \frac{-1 + 12z - 24z^2 + (1 - 8z)^{3/2}}{32z^2}.
\end{equation}
Then the estimate for the number of bipartite maps with $n$ edges is
\begin{equation*}
	[z^n]A_b(z) \sim \frac{3}{2\sqrt{\pi}}\, n^{-5/2} \, 8^n
	\qquad\text{as } n\to\infty.
\end{equation*}

The next step is also rather direct since a map (or graph) is bipartite if and only if its maximal 2-connected components are bipartite.
The decomposition of a map into 2-connected components implies the following equation for the series $B_b$ of 2-connected bipartite maps:
\begin{equation*}
	A_b(z) = 1 + zA_b(z)^2 + B_b(zA_b(z)^2).
\end{equation*}
Notice that the single loop is not a bipartite map.
By elimination one obtains an algebraic equation for $B_b(z)$ of degree 5 from which (see \cite{LW04}) one deduces the estimate
\begin{equation*}
	[z^n]B_b(z)\sim \frac{75\sqrt{65}}{4394\sqrt{\pi}}\, n^{-5/2}\, \left(\frac{128}{25}\right)^n
	\qquad\text{as } n\to\infty.
\end{equation*}
However, the next step, i.e. to encode 3-connected maps starting from  2-connected maps, does not generalise at all since it is far from true that a 2-connected graph is bipartite if and only if its 3-connected components are bipartite (consider for instance the graph obtained from $K_4$ by subdividing once every edge).

Our approach is to enrich the framework by counting \textit{bicoloured} maps, which are maps together with a 2-colouring (black and white) of the vertices such that the root vertex is always black.
The 2-colouring is arbitrary and not necessarily proper.
Let $\mathcal{C}$ be a class of bicoloured maps with at least one edge and let $C(z,\nu)$ be the associated series, where $z$ marks all edges as before and $\nu$ marks monochromatic edges, that is, edges whose endpoints have the same colour.
We will consider the decomposition $C(z,\nu) = C_1(z,\nu) + C_2(z,\nu)$, where $C_1(z,\nu)$ counts the maps in $\mathcal{C}$ rooted at a monochromatic edge and $C_2(z,\nu)$ counts those rooted at a bichromatic edge.
It is then clear that $C_2(z,0)$ is the series of bipartite maps in $\mathcal{C}$.
Our strategy is to compute $C_2(z,\nu)$ for general, 2-connected and 3-connected maps.

After his work on map enumeration in the 1960's, Tutte went back to the original motivation of  counting 4-colourings of triangulations.
We recall that the \textit{chromatic polynomial} $\chi_{\mathsf{g}}(q)$ of a graph $\mathsf{g}$ is a function equal to the number of proper $q$-colourings of $\mathsf{g}$.
For a class of maps $\mathcal{C}$ Tutte defined its `chromatic sum' as
\begin{equation}\label{eq:chromatic_sum}
	C(z,q)=	\sum_{\mathsf{m}\in \mathcal{C}} z^{|\mathsf{m}|} \chi_{\mathsf{m}}(q),
\end{equation}
where $|\mathsf{m}|$ denotes the number of edges of $\mathsf{m}$.
If one could find a way of computing $C(z,q)$, then one could compute the average number of $q$-colourings of maps in $\mathcal{C}$.
Tutte devoted a long series of papers to chromatic sums of triangulations for different values of $q$.
Starting with~\cite{T73}, Tutte obtained that for any $q$, the generating function for 3-connected triangulations weighted by their chromatic polynomial satisfies a differential equation from which one can obtain explicit recurrence relations (see~\cite{T82, odl83}).
For $q=4$, the differential equation obtained is somehow simple (see Equation (52) in~\cite{T82}), so recursive enumerative formulas exist in this case. However, as he remarked years later in~\cite{T98}
\begin{quote}
	I said near the beginning of this chapter that information about averages might be easier to obtain than a proof of the Four Colour Theorem.
	Yet now we have the Haken-Appel proof, and we still lack an explicit formula or even an asymptotic approximation for our four-colour average.
\end{quote}
In fact, nowadays to find the asymptotic number of $4$-coloured triangulations is still an open problem.
The topic was revived much later by Bernardi and Bousquet-M\'elou in a remarkable paper~\cite{BBM11} (and its sequel~\cite{BBM17}), where they write:
\begin{quote}
	This tour de force has remained isolated since then, and it is our objective to reach a better understanding of Tutte's rather formidable approach, and to apply it to other problems in the enumeration of colored planar maps.
\end{quote}

The approach in~\cite{BBM11} is to consider all $q$-colourings, proper or not, and take as a new parameter the number of monochromatic edges.
Our starting point is Theorem 21 from~\cite{BBM11}, where the authors determine the series $M(z,\nu)$ of bicoloured maps, where $z$ marks edges and $\nu$ monochromatic edges.
It is an algebraic function of degree 6 whose first terms are
\begin{equation*}
	M(z,\nu) = 1 + \left(1 + 2\,\nu \right)z + \left(9\,\nu^2 + 8\,\nu + 3\right)z^2 + \left(42\,\nu^3 + 72\,\nu^2 + 51\,\nu + 12\right)z^3 + \cdots
\end{equation*}
For instance, the term $z$ corresponds to the isthmus map with the unique proper 2-colouring (recall that the root vertex is by convention coloured black) and $2\nu z$ corresponds to the loop and isthmus maps coloured with a single color.

In Section~\ref{sec:enum_bicolor_maps} we obtain the series $M_1(z,\nu)$ of bicoloured maps where the root edge is monochromatic, and as a consequence the series $M_2(z,\nu) = M(z,\nu) - 1 - M_2(z,\nu)$ of bicoloured maps rooted at a bichromatic edge.
Once $M_1$ and $M_2$ are determined, using the decomposition of maps into 2-connected components we obtain the generating functions $B_1(z,\nu)$ and $B_2(z,\nu)$ of 2-connected bicoloured maps.
Then using the decomposition of 2-connected maps into 3-connected components we obtain the generating functions $T_1(z,\nu)$ and $T_2(z,\nu)$ of 3-connected bicoloured maps.

\begin{figure}
	\begin{center}
		\begin{minipage}{.3\textwidth}
			\centering
			\small{\textsf{12 edges}}
			\includegraphics[scale=1]{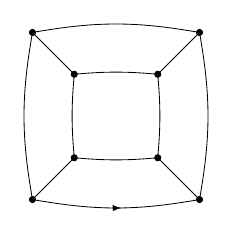}\\
			\vspace{1cm}\small{\textsf{18 edges}}
			\includegraphics[scale=1]{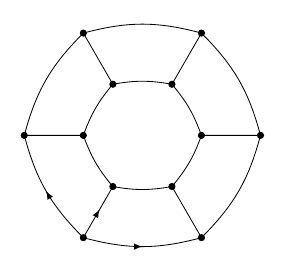}\\
			\vspace{1cm}\small{\textsf{20 edges}}
			\includegraphics[scale=1]{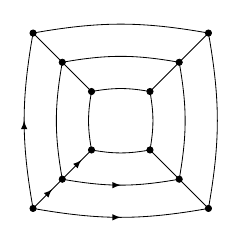}
		\end{minipage}
		\begin{minipage}{.3\textwidth}
			\centering
			\vspace{-.2cm}
			\small{\textsf{16 edges}}
			\includegraphics[scale=1]{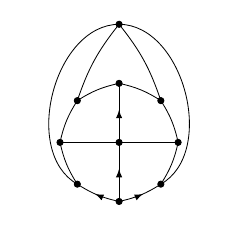}\\
			\vspace{1.2cm}\small{\textsf{18 edges}}
			\includegraphics[scale=1]{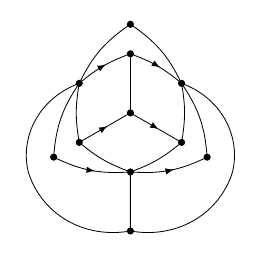}\\
			\vspace{1.2cm}\small{\textsf{20 edges}}
			\includegraphics[scale=1]{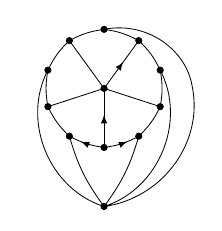}
		\end{minipage}
		\begin{minipage}{.34\textwidth}
			\centering
			\vspace{5.2cm}
			\small{\textsf{19 edges}}\\
			\includegraphics[scale=1]{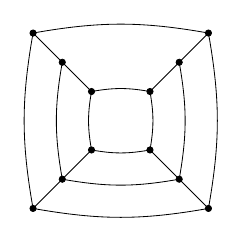}\\
			\vspace{1.2cm}\vspace{.3cm}
			\small{\textsf{20 edges}}
			\includegraphics[scale=1]{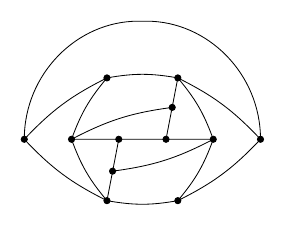}
		\end{minipage}
	\end{center}
	\caption{
		The smallest 3-connected bipartite planar maps pictured with their rootings and number of edges.
		The cube is the unique (unrooted) 3-connected bipartite map with 12 edges; it has 24 symmetries and thus admits a unique rooting.
		There are no 3-connected bipartite maps with 13, 14 and 15 edges.
		The  unique 3-connected bipartite map with 16 edges has 8 symmetries, hence $[z^{16}]T_b(z) = 2\cdot 16/8 = 4$.
		There are two 3-connected bipartite maps with 18 edges, the hexagonal prism (left) has 12 symmetries and the other one has 6.
		The unique 3-connected bipartite map with 19 edges has only 2 symmetries.
		The three 3-connected bipartite maps with 20 edges admit 8 (left), 10 (middle) and 2 (right) symmetries.
		The rootings of the two maps with only 2 symmetries (middle right and bottom right) are not drawn for readability.
	}\label{fig:3-conn_bipartite}
\end{figure}

Finally, $T_b(z) = T_2(z,0)$ is the series of 3-connected bipartite maps, which is what we needed.
$T_b(z)$ is algebraic of degree 26 and its minimal polynomial is too long to be reproduced here, but can be found, together with the relevant computations regarding this work,  in the accompanying \texttt{Maple} session~\cite{MapleSession}.
It is somehow surprising that the series of such a natural class of planar maps has this algebraic complexity.
The first terms are
\begin{equation*}
	T_b(z) = z^{12} + 4z^{16} + 9z^{18} + 19z^{19} + 29z^{20} + 63z^{21} +  198z^{22} +  345z^{23} +  685z^{24} +1775z^{25} +  \cdots,
\end{equation*}
Recall that the coefficient of $z^n$ counts the different rootings of 3-connected bipartite unrooted maps with $n$ edges, where the number of rootings of a map is twice the number of edges divided by the number of symmetries.
In particular, let us verify the first five coefficients following the order left to right then top to bottom of their illustration in Figure~\ref{fig:3-conn_bipartite}:
Remark that the hexagonal prism is the smallest 3-connected bipartite planar map that is not a quadrangulation (compare with the last table in~\cite{NRR19}).

The computations are within the capabilities of a modern computer algebra system, and we can perform the corresponding singularity analysis to obtain our main result.
\begin{theorem}\label{th:main}
	Let $t_n$ be the number of 3-connected bipartite maps with $n$ edges.
	Then
	\begin{equation*}
		t_n\sim t\, n^{-5/2}\, \gamma^n
		\qquad\text{as } n\to\infty,
	\end{equation*}
with $t \approx 0.00412$ and $\gamma = \rho^{-1} \approx 2.40958$, and where $\rho \approx 0.41501$ is a root of the irreducible polynomial
	\begin{equation*}
		\begin{array}{ll}
			82796609536z^{10} - 125942890496z^9 - 278408523776z^8 + 430685329920z^7 + 381960569664z^6 \\
			- 519574865472z^5 - 253658186064z^4 + 253489557672z^3 + 78314553945z^2 \\
			- 51532664454z + 1113912729.
		\end{array}
	\end{equation*}
\end{theorem}
The class of 3-connected bipartite planar maps is a naturally defined class of planar maps that has not been counted before.
Other natural classes of maps for which the enumeration problem remains open include those that are 4-connected, triangle-free, or 3-colourable.

We conclude this introduction by showing the growth constants of some of the classes of maps discussed above, in all cases counted by number of edges.
The fact that the first two values in the third row are the same is because a connected bipartite cubic graph is necessarily 2-connected.
\begin{center}
{
	\renewcommand{\arraystretch}{1.3}
	\begin{tabular}{l|ccc}
		Class of maps & Arbitrary & 2-connected & 3-connected \\
		\hline
		Arbitrary & 12 & $27/4$ & 4 \\
		Bipartite & 8  & $125/8$ & $\gamma \approx 2.40958$ \\
		Bipartite cubic & 2& 2 & $8/5$
	\end{tabular}
}
\end{center}
%

%%%%%%%%%%%%%%%%%%%%%%%%%%%%%%%%%%%%%%%%%%%%%%%%%%%%%%%%%%%%%%%%%%%%%%%%%%%%%%%%
%
%	Counting bicoloured 3-connected planar maps
%
%%%%%%%%%%%%%%%%%%%%%%%%%%%%%%%%%%%%%%%%%%%%%%%%%%%%%%%%%%%%%%%%%%%%%%%%%%%%%%%%
\section{Counting bicoloured 3-connected planar maps}\label{sec:enum_bicolor_maps}

In this section we first determine the generating function of bicoloured maps rooted at a monochromatic edge or at a bichromatic edge.
Then, using the decomposition of maps into maximal 2-and 3-connected components, we determine the series of 2-connected and 3-connected bicoloured maps.

\subsection{Counting bicoloured maps}

\paragraph{The Ising polynomial of a map.}

The partition function of the Potts model (a model in statistical physics for spin interactions) on a graph $\mathsf{m}$ is defined as
\begin{equation*}
	P_{\mathsf{g}}(q,\nu) = \sum_{c\colon V(\mathsf{g})\to [q]} \nu^{m(c)},
\end{equation*}
where $[q] = \{1,2,\dots,q\}$, the sum is defined over all vertex colourings (with $q$ colours) of $\mathsf{g}$, and $m(c)$ is the number of monochromatic edges defined by $c$.
When restricted to two colours, $P_{\mathsf{g}}(2,\nu)$ is also known as the partition function of the Ising model (a model for ferromagnetism, where the two colours correspond to the two possible  values of the ``spin'' of a particle).

The function $P_{\mathsf{g}}(q,\nu)$ has some useful properties.
Let $\mathsf{g}_1$ and $\mathsf{g}_2$ be two vertex-disjoint graphs, and let $\mathsf{g}$ be the graph obtained by identifying a vertex of $\mathsf{g}_1$ and a vertex of $\mathsf{g}_2$.
Then
\begin{equation}\label{potts_properties}
	P_{\mathsf{g}_1\cup \mathsf{g}_2}(q,\nu) = P_{\mathsf{g}_1}(q,\nu)P_{\mathsf{g}_2}(q,\nu)
	\quad\text{and}\quad
	P_{\mathsf{g}}(q,\nu) = \frac{1}{q}P_{\mathsf{g}_1}(q,\nu)P_{\mathsf{g}_2}(q,\nu).
\end{equation}
Furthermore, $P_{\mathsf{g}}(q,\nu)$ can be computed recursively as follows.
If $\mathsf{g}$ is the empty graph then $P_{\mathsf{g}}(q,\nu) = q^{|V(\mathsf{g})|}$.
Otherwise we have
\begin{equation}\label{eq:induction_potts}
	P_{\mathsf{g}}(q,\nu) = P_{\mathsf{g}\backslash e}(q,\nu) + (\nu - 1)P_{\mathsf{g}/e}(q,\nu),
\end{equation}
where $\mathsf{g}\backslash e$ is the graph obtained by deleting the edge $e$ from $\mathsf{g}$ and $\mathsf{g}/e$ is obtained by contracting $e$ (when $e$ is a loop, $\mathsf{g}/e$ is obtained by deleting $e$).
As pointed out in~\cite[Section 2.3]{BM11}, a combinatorial explanation of this equation is that $\nu P_{\mathsf{g}/e}(q,\nu)$ counts colourings of $\mathsf{g}$ for which the edge $e$ is monochromatic while $P_{\mathsf{g}\backslash e}(q,\nu) - P_{\mathsf{g}/e}(q,\nu)$ counts those where $e$ is bichromatic, as we do not want to contract a bichromatic edge.
Observe that the relation~\eqref{eq:induction_potts} implies by induction that $P_{\mathsf{g}}(q,\nu)$ is a polynomial in $q$ and $\nu$ with no constant term, in particular it is divisible by $q$.
Observe then that $P_{\mathsf{g}}(q,\nu)$ encodes all $q$-colourings of $\mathsf{g}$, while $P_{\mathsf{g}}(q,\nu)/q$ will encode those whose root vertex is of a fixed colour.
$P_{\mathsf{g}}(q,\nu)$ is in fact equivalent to the ubiquitous Tutte polynomial $T_{\mathsf{g}}(x,y)$ through the simple change of variables $y=\nu$ and $(x - 1)(y - 1) = q$,
while $P_{\mathsf{g}}(q,0)$ is the chromatic polynomial of $\mathsf{g}$ (see~\cite{BM11} for an extended account).
We call $P_{\mathsf{g}}(2,\nu)$ the Ising polynomial $\mathsf{g}$.

\begin{figure}[h]
	\begin{center}
	\begin{minipage}{.58\textwidth}
	\centering
		\includegraphics[scale=1]{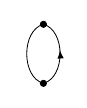}
			\raisebox{.75cm}{$\overset{\textsf{m}_/}{\boldsymbol{\dashleftarrow}}$}
		\includegraphics[scale=1]{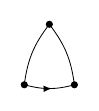}
			\raisebox{.75cm}{$\overset{\textsf{m}_{\backslash}}{\boldsymbol{\dashrightarrow}}$}
		\includegraphics[scale=1]{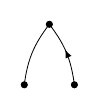}\\
		\includegraphics[scale=1]{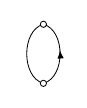}
			\raisebox{.75cm}{$\overset{\textsf{m}_/}{\boldsymbol{\dashleftarrow}}$}
		\includegraphics[scale=1]{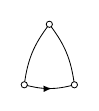}
			\raisebox{.75cm}{$\overset{\textsf{m}_{\backslash}}{\boldsymbol{\dashrightarrow}}$}
		\includegraphics[scale=1]{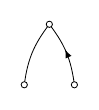}\\
		\includegraphics[scale=1]{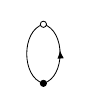}
			\raisebox{.75cm}{$\overset{\textsf{m}_/}{\boldsymbol{\dashleftarrow}}$}
		\includegraphics[scale=1]{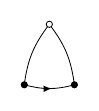}
			\raisebox{.75cm}{$\overset{\textsf{m}_{\backslash}}{\boldsymbol{\dashrightarrow}}$}
		\includegraphics[scale=1]{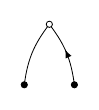}\\
		\includegraphics[scale=1]{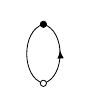}
			\raisebox{.75cm}{$\overset{\textsf{m}_/}{\boldsymbol{\dashleftarrow}}$}
		\includegraphics[scale=1]{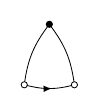}
			\raisebox{.75cm}{$\overset{\textsf{m}_{\backslash}}{\boldsymbol{\dashrightarrow}}$}
		\includegraphics[scale=1]{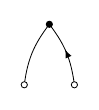}
	\end{minipage}
	\hfill\vline\hfill
	\begin{minipage}{.38\textwidth}
	\centering
		\includegraphics[scale=1]{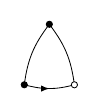}
			\raisebox{.75cm}{$\overset{\textsf{m}_{\backslash}}{\boldsymbol{\dashrightarrow}}$}
		\includegraphics[scale=1]{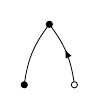}\\
		\includegraphics[scale=1]{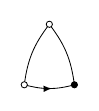}
			\raisebox{.75cm}{$\overset{\textsf{m}_{\backslash}}{\boldsymbol{\dashrightarrow}}$}
		\includegraphics[scale=1]{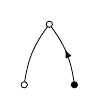}\\
		\includegraphics[scale=1]{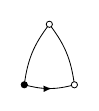}
			\raisebox{.75cm}{$\overset{\textsf{m}_{\backslash}}{\boldsymbol{\dashrightarrow}}$}
		\includegraphics[scale=1]{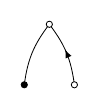}\\
		\includegraphics[scale=1]{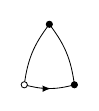}
			\raisebox{.75cm}{$\overset{\textsf{m}_{\backslash}}{\boldsymbol{\dashrightarrow}}$}
		\includegraphics[scale=1]{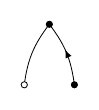}
	\end{minipage}
	\end{center}
	\caption{
		The eight different 2-colourings of the rooted triangle map $\textsf{m}$, that are encoded by the Ising polynomial $P_{\textsf{m}}(2,\nu)$, and their decompositions.
		On the left are the 2-colourings of $\textsf{m}$ for which the root edge is monochromatic.
		Their contribution to $P_{\textsf{m}}(2,\nu)$ is $2\nu^3 + 2\nu$.
		Alongside them are the 2-coloured maps obtained after contracting (left) or deleting (right) their root edge.
		On the right are the 2-colourings of $\textsf{m}$ for which the root edge is bichromatic.
		Their contribution to $P_{\textsf{m}}(2,\nu)$ is $4\nu$.
		Alongside them are the 2-coloured maps obtained after deleting their root edge (recall that we do not contract a bichromatic root edge).
		Thus only the left side contributes to $P_{\textsf{m}_/}(2,\nu) = 2\nu^2 + 2$, while both sides contribute to $P_{\textsf{m}_{\backslash}}(2,\nu) = 2\nu^2 + 4\nu + 2$.
		Observe that we indeed have $P_{\textsf{m}}(2,\nu) = P_{\textsf{m}_{\backslash}}(2,\nu) + (\nu - 1)P_{\textsf{m}_/}(2,\nu) =  2\nu^3 + 6\nu$.
	}\label{fig:ising_maps}
\end{figure}

Additionally, \eqref{eq:induction_potts} allows us to define by induction the \textit{Ising polynomial $P_{\textsf{m}}(2,\nu)$ of a rooted map $\textsf{m}$}, as follows.
The Ising polynomial of the empty map is equal to 2.
The deletion and contraction operations are performed on the root edge $e$, and the resulting maps are denoted by $\textsf{m}_{\backslash}$ and $\textsf{m}_/$, respectively.
The new root edge is defined canonically as the next edge encountered while walking along the boundary of the root face of $\textsf{m}$ in the direction induced by $e$.
In the case where $e$ is a bridge of $\textsf{m}$, then $\textsf{m}_{\backslash}$ is  composed of two disjoint maps $\textsf{n}$ and $\textsf{n}'$ and $P_{\textsf{m}_{\backslash}}(2,\nu) = P_{\textsf{n}}(2,\nu)P_{\textsf{n}'}(2,\nu)$.
An illustration is given in Figure~\ref{fig:ising_maps}.

\paragraph{The Ising generating function of maps.}

In this section we use the common terminology in~\cite{BBM11} and~\cite{BM11}.
Let $\mathcal{M}$ be the class of all planar maps, and let $\mathcal{M}_0$ be equal to $\mathcal{M}$ minus the empty map.
The generating function $M(z,\nu)$ of bicoloured maps can be obtained directly by weighting each map in $\mathcal{M}$ by its Ising polynomial (compare with~\eqref{eq:chromatic_sum}).
Then we have
\begin{equation*}
	M(z,\nu) = \frac{1}{2}\sum_{\mathsf{m}\in\mathcal{M}} P_\mathsf{m}(2,\nu)\, z^{|\mathsf{m}|},
\end{equation*}
where  $|\mathsf{m}|$ is the number of edges in $\mathsf{m}$, and the factor $1/2$ is because the root vertex is coloured black.
The relation~\eqref{eq:induction_potts} applied to $P_\mathsf{m}(2,\nu)$ induces a functional equation for $M(z,\nu)$ analogous to the one for (uncoloured) maps in~\eqref{eq:induction_maps}.
In the same way that for edge deletion in~\eqref{eq:induction_maps} we encoded the degree of the root face in the  variable $y$, the operation of edge contraction in~\eqref{eq:induction_potts} requires a {second}  variable $x$ encoding the degree of the root vertex.
We remark that  $x$ and $y$ play dual r\^oles: the degree of the root face of a map $\mathsf{m}$ is the degree of the root vertex of the geometric dual $\mathsf{m}^*$.
Also, $(\mathsf{m} \backslash e)^* \simeq \mathsf{m}^* / e^*$, where $e^*$ is the dual edge of $e$, and $(\mathsf{m} / e)^* \simeq \mathsf{m}^* \backslash e^*$.
In terms of the Tutte polynomial, this implies that $T_{\mathsf{m}^*}(x,y) = T_\mathsf{m}(y,x)$.

We write now $N(x,y) = M(x,y; z,\nu)$ in order to emphasize the role of the catalytic variables.
Let  $\mathrm{dv}(\mathsf{m})$ and $\mathrm{df}(\mathsf{m})$ be the degree of the root vertex and the degree of the root face of $\mathsf{m}$, respectively.
Following~\cite{BM11}, we define
\begin{equation*}
	N_{\backslash}(x,y)= \frac{1}{2}\sum_{\mathsf{m}\in \mathcal{M}_0} z^{|\mathsf{m}|} x^{\mathrm{dv}(\mathsf{m})} y^{\mathrm{df}(\mathsf{m})} P_{\mathsf{m}_{\backslash}}(2,\nu)
	\quad\text{and}\quad
	N_{/}(x,y)= \frac{1}{2}\sum_{\mathsf{m}\in \mathcal{M}_0} z^{|\mathsf{m}|} x^{\mathrm{dv}(\mathsf{m})} y^{\mathrm{df}(\mathsf{m})} P_{\mathsf{m}_/}(2,\nu).
\end{equation*}
The relation~\eqref{eq:induction_potts} directly implies that
\begin{equation*}
	N(x,y) = 1 + N_{\backslash}(x,y) + (\nu - 1)N_{/}(x,y).
\end{equation*}
Furthermore, generalising~\eqref{eq:induction_maps} with the help of~\eqref{potts_properties}, one can write $N_{\backslash}(x,y)$ and $N_{/}(x,y)$ as functions of $N(x,y)$ and its evaluations (see~\cite[Proposition 5.1]{BM11} for a detailed proof):
\begin{equation}\label{eq:indutcion_maps_ising_deletion}
	N_{\backslash}(x,y) = 2xy^2zN(1,y)N(x,y) + (x - 1)xyzN(x,1)N(x,y) + xyz\frac{yN(x,y) - N(x,1)}{y - 1}
\end{equation}
and
\begin{equation}\label{eq:indutcion_maps_ising_contraction}
	N_{/}(x,y) = x^2yzN(x,1)N(x,y) + (y - 1)xyzN(1,y)N(x,y) + xyz\frac{xN(x,y) - N(1,y)}{x - 1}.
\end{equation}
Combining~\eqref{eq:indutcion_maps_ising_deletion} and~\eqref{eq:indutcion_maps_ising_contraction} gives
\begin{equation}\label{eq:induction_maps_ising}
	\begin{split}
		N(x,y) = 1
		& + xyz(2y + (\nu - 1)(y - 1))N(x,y)N(1,y) + xyz(x\nu - 1)N(x,y)N(x,1) \\
		& + xyz(\nu - 1)\frac{xN(x,y) - N(1,y)}{x - 1} + xyz\frac{yN(x,y) - N(x,1)}{y - 1}.
	\end{split}
\end{equation}

Equation~\eqref{eq:induction_maps_ising} was first established in an equivalent form for maps weighted by their Tutte polynomial (called the \textit{dichromate} by Tutte) in~\cite{T71}.
Using a remarkable method to solve this particular family (for general values of $q$) of functional equations with two catalytic variables, it is proven in~\cite{BBM11} that $M(z,\nu)$ is an algebraic function of degree 6 which admits the rational parametrisation
\begin{equation}\label{eq:param_M}
	\begin{split}
		M(z,\nu) = \frac{1 + 3\nu S - 3\nu S^2 - \nu^2S^3}{(1 - 2S + 2\nu^2S^3 -\nu^2S^4)^2}\Big(1
		& - (3 + \nu)S + (1 + 2\nu)S^2 - \nu(1 - 5\nu)S^3 \\
		& + \nu(1 - 6\nu)S^4 + 2\nu^2(1 - \nu)S^5 + \nu^3S^6\Big),
	\end{split}
\end{equation}
where $S = S(z,\nu)$ is defined as the unique solution of
\begin{equation}\label{eq:param_S}
	S = z\frac{(1 + 3\nu S - 3\nu S^2 - \nu^2 S^3)^2}{1 - 2S + 2\nu^2S^3 - \nu^2S^4},
\end{equation}
whose coefficients at $z=0$ are polynomials in $\nu$ with non-negative coefficients.

\paragraph{Rooting at a monochromatic edge.}

Let us denote by $M_1(z,\nu)$ and $M_2(z,\nu)$ the generating functions for bicoloured maps where the root edge is monochromatic and bichromatic, respectively, and where $z$ marks edges and $\nu$ monochromatic edges, so that $M(z,\nu) = 1 + M_1(z,\nu) + M_2(z,\nu)$.
As observed above, $\nu P_{\textsf{m}_/}(z,\nu)$ encodes colourings of $\textsf{m}$ for which the root edge is monochromatic while $P_{\textsf{m}_{\backslash}}(z,\nu) - P_{\textsf{m}_/}(z,\nu)$ encodes those where the root edge is bichromatic.
Hence
\begin{equation*}	
	M_1 (z,\nu) = \nu  N_/(1,1)
	\quad\text{and}\quad
  	M_2(z,\nu)=  N_{\backslash} (1,1) - N_/(1,1).
\end{equation*}
By setting $y=1$ in~\eqref{eq:indutcion_maps_ising_contraction} we get
\begin{equation*}
	N_/(x,1) = x^2zN(x,1)^2 + xz\frac{xN(x,1) - N(1,1)}{x - 1}.
\end{equation*}
Note that the limit as $x\to 1$ of the last summand of the right hand side of the above equation is the derivative of $xN(x,1)$ evaluated at $x=1$.
So that taking the limit as $x\to 1$ of the above equation gives
\begin{align*}
	\underset{x\to 1}{\lim} \,\, N_/(x,1) = zN(1,1)^2 + zN(1,1) + zN_x(1,1),
	\quad\text{where } N_a(1,1) = \frac{\partial N(x,y)}{\partial a}{\Bigr|}_{x=1,y=1}
	\text{ for } a\in\{x,y\},
\end{align*}
which can be rephrased using the identities $M_1 (z,\nu) = \nu  N_/(1,1)$ and $M(z,\nu) = N(1,1)$ into
\begin{equation}\label{eq:M1}
	M_1(z,\nu) = \nu zM(z,\nu)^2 + \nu zM(z,\nu) + \nu zN_x(1,1).
\end{equation}

In~\cite{BBM11}, the authors developed a method to first eliminate the variable $x$ in Equation~\eqref{eq:induction_maps_ising} then the variable $y$, resulting in two polynomials equations $R_1(X,Y,z,\nu)$ and $R_2(X,Y,z,\nu)$ such that $R_1=R_2=0$ when $X = M(z,\nu)$ and $Y = N_y(1,1)$ (see in particular the proof of~\cite[Theorem 21]{BBM11}).
From there, eliminating $N_y(1,1)$ leads to the minimal polynomial of $M(z,\nu)$.
By a direct step by step adaptation of the proof of~\cite[Theorem 21]{BBM11}, but instead eliminating the variable $y$ first then the variable $x$ (see the \texttt{Maple} session~\cite{MapleSession}), we obtain a polynomial equation $R(X,Y,z,\nu)$ such that $R=0$ when $X = M(z,\nu)$ and $Y = N_x(1,1)$ (in fact we obtain two, but one is enough for our purpose).

Finally, eliminating from the system composed of Equations~\eqref{eq:param_M},~\eqref{eq:param_S},~\eqref{eq:M1} and $R=0$ yields an irreducible polynomial equation $Q^{(1)}(M_1,z,\nu) = 0$, also of degree six in $M_1$, defining $M_1 = M_1(z,\nu)$ implicitly as a function of $z$ and $\nu$.
An equation $Q^{(2)}(M_2,z,\nu) = 0$ can be derived by elimination, using $M(z,\nu) = 1 + M_1(z,\nu) + M_2(z,\nu)$.
The polynomial $Q^{(2)}$ has also degree six in $M_2 = M_2(z,\nu)$.
Incidentally, both $M_1(z,\nu)$ and $M_2(z,\nu)$ also admit, as $M(z,\nu)$, rational parametrisations in $S$ given by
\begin{align*}
	& M_1(z,\nu) = \frac{\nu S}{(1 - 2S + 2\nu^2S^3 -\nu^2S^4)^2}
	\Big(\nu^4S^8 + \nu^3(5 - 2\nu)S^7 + \nu^2(7 - 16\nu)S^6 + \nu(14\nu^2 - 24\nu + 3)S^5 \\
	& \qquad\qquad\qquad\qquad\qquad\qquad\qquad  + \nu(33\nu - 8)S^4 - (16\nu^2 - 7\nu - 2)S^3 - (8\nu + 7)S^2 + (3\nu + 8)S - 2\Big), \\
	& M_2(z,\nu) = \frac{S}{(1 - 2S + 2\nu^2S^3 -\nu^2S^4)^2}
	\Big(\nu^4S^7 - 3\nu^4S^6 + \nu^3(\nu - 2)S^5 + \nu^2(4\nu + 5)S^4 - 5\nu^2S^3 \\
	& \qquad\qquad\qquad\qquad\qquad\qquad\qquad - \nu(\nu + 4)S^2 + (2\nu + 3)S - 1\Big).
\end{align*}

\subsection{Counting bicoloured 2- and 3-connected maps}

Let $B_1(z,\nu)$ and $B_2(z,\nu)$ be the series counting 2-connected bicoloured maps rooted at a monochromatic and bichromatic edge, respectively.
Then a straightforward generalisation of~\eqref{eq:MN} gives
\begin{align}
	M_1(z,\nu) & = B_1(zM(z,\nu)^2,\nu) + 2\nu z M(z,\nu)^2, \label{eq:B1} \\
	M_2(z,\nu) & = B_2(zM(z,\nu)^2,\nu) +  zM(z,\nu)^2. \label{eq:B2}
\end{align}
Those equations reflect the fact that the root edge in the 2-core $\mathsf{C}(\mathsf{m})$ of a map $\mathsf{m}$ is actually the root edge of $\mathsf{m}$, hence $\mathsf{m}$ and its 2-core are either both monochromatic or both bichromatic.
The substitutions  are in terms of $M(z,\nu)$ only, since the maps placed at the corners of $\mathsf{C}(\mathsf{m})$ share no edge with $\mathsf{C}(\mathsf{m})$.
By elimination from Equations~\eqref{eq:param_M},~\eqref{eq:param_S}, $Q^{(1)}=0$ and~\eqref{eq:B1}, we obtain the minimal polynomial of $B_1(z,\nu)$ of degree 13.
The same is done from Equations~\eqref{eq:param_M},~\eqref{eq:param_S}, $Q^{(2)}=0$ and~\eqref{eq:B2} to obtain the minimal polynomial of $B_2(z,\nu)$, also of degree 13.
Both polynomials can be found in~\cite{MapleSession}.

Finally, we determine the series $T_2(z,\nu)$ counting 3-connected bicoloured maps rooted at a bichromatic edge, and incidentally the series $T_1(z,\nu)$ counting those rooted at a monochromatic edge.
Let $D_1 = B_1(z,\nu)/(z\nu)$ and $D_2 = B_2(z,\nu)/z$.
Let $S_i = S_i(z,\nu)$ and $P_i = P_i(z,\nu)$ encode series and parallel compositions as for bicoloured maps, where the index $i$ has the same meaning as before.
Then a straightforward generalisation of~\eqref{eq:BT} gives
\begin{equation}\label{eq:D1D2}
	D_1 = z\nu + S_1 + P_1 + T_1\left(D_2, \frac{D_1}{D_2}\right),
	\qquad\qquad
	D_2 = z + S_2 + P_2 + T_2\left(D_2, \frac{D_1}{D_2}\right).
\end{equation}
The terms $T_i\left( D_2, D_1/D_2 \right)$ encode the replacement by 2-connected maps of the non-root edges of the 3-cores $\mathsf{T}(\mathsf{m})$, where monochromatic (resp. bichromatic) edges of $\mathsf{T}(\mathsf{m})$ have to be replaced by maps rooted at a monochromatic (resp. bichromatic) edge.
We also have the following relations:
\begin{equation} \label{eq:PiSi}
	\renewcommand{\arraystretch}{1.6}
	\begin{array}{llllll}
		S_1 = (D_1 - S_1)D_1 + (D_2 - S_2)D_2, & \qquad
		S_2 = (D_2 - S_2)D_1 + (D_1 - S_1)D_2, \\
		P_1 = (D_1-P_1)D_1 & \qquad P_2 = (D_2-P_2)D_2. \\
	\end{array}
\end{equation}
For the equation for $S_1$, remark that in order to obtain a series map rooted at a monochromatic edge one must compose in series two maps that are either both rooted at a monochromatic edge or both rooted at a bichromatic edge.
While for $S_2$ one root must be monochromatic and the other bichromatic.
The equations for $P_1$ and $P_2$ are simpler since both root edges in the parallel compositions must be of the same kind.

Eliminating from the system composed of~\eqref{eq:D1D2},~\eqref{eq:PiSi} and the minimal polynomials of $B_1(z,\nu)$ and $B_2(z,\nu)$, one obtains the minimal polynomial of $T_2 = T_2(z,\nu)$, which turns out to be of degree 26 in $T_2$ (see~\cite{MapleSession}).

%%%%%%%%%%%%%%%%%%%%%%%%%%%%%%%%%%%%%%%%%%%%%%%%%%%%%%%%%%%%%%%%%%%%%%%%%%%%%%%%
%
%	Proof of Theorem 1
%
%%%%%%%%%%%%%%%%%%%%%%%%%%%%%%%%%%%%%%%%%%%%%%%%%%%%%%%%%%%%%%%%%%%%%%%%%%%%%%%%
\section{Proof of Theorem \ref{th:main}}

We have obtained $T_2(z,\nu)$ in the previous section and, as discussed in the introduction, the series of 3-connected bipartite maps is equal to $T_b(z) = T_2(z,0)$.
The minimal polynomial of $T_b(z)$ is of the form
\begin{equation*}
	P(z,T_b) = \sum_{i=0}^{26} p_i(z) T_b^i,
\end{equation*}
where the $p_i(z)$'s are polynomials in $z$ which happen to be of degree $85-i$.
In particular, the leading coefficient of $P(z,T_b)$ is given by
\begin{equation}\label{eq:leading_coeff_T}
	p_{26}(z) = 4096(z + 3)^{11}(z - 1)^{22}(z + 1)^{26}.
\end{equation}
Because it is algebraic, $T_b(z)$ can be represented at $z=0$, i.e. in a neigbourhood of $0$, as a generating function with non-negative coefficients and radius of convergence $\rho$, for some $\rho > 0$, corresponding to a branch of the curve $P(z,T_b) = 0$ passing through the origin.
We refer the reader to the detailed discussion in~\cite[Chapter VII.7]{FS09} for this and related facts on algebraic generating functions used in the rest of the proof.

Next we find the value of $\rho$.
Note first that since the radius of convergence of the series of all 3-connected maps is $1/4$, as discussed in Section~\ref{sec:introduction}, we must have $\rho \ge 1/4$.
In order to get an upper bound for $\rho$ we need a subclass of 3-connected bipartite maps which is large enough.
This is provided by the class of 3-connected bipartite cubic maps (called bicubic maps by Tutte in~\cite{T63}).
As shown in~\cite{MRS74} the radius of convergence of this class when maps are counted by the number of edges is $5/8$.
It follows that $\rho\in (1/4,5/8)$.
Second, by Pringsheim's theorem (see~\cite[Theorem IV.6]{FS09}) $\rho$ must be a singularity of $T_b(z)$.
And since $T_b(z)$ is algebraic its singularities can be of two types: they are either \textit{algebraic poles}, i.e. points for which the degree of $P(z,T_b)$ decreases (by cancelling its leading coefficient), or they are \textit{branch points}, that are roots of the discriminant $d(z)$ of $P(z,T_b)$ with respect to $T_b$.
It is clear that $\rho$ has to be a branch point singularity as none of the roots of~\eqref{eq:leading_coeff_T} are in the interval $(1/4,5/8)$.
Let us now consider the polynomial $d(z)$.
It has 242 different roots but only one lays in the interval $(1/4,5/8)$: the root approximately equal to $0.41501$.
It must then be the radius of convergence of $T_b(z)$, and the unique irreducible factor of $d(z)$ having $\rho\approx 0.41501$ as a root is precisely the one claimed in Theorem~\ref{th:main}.

In order to obtain an estimate for $[z^n]T_b(z)$ we  apply the following \textit{transfer theorem}~\cite[Chapter VI.3]{FS09}.
Assume that $F(z)$ has radius of convergence $\phi > 0$ and is analytic in an open domain at $z = \phi$ of the form
\begin{equation*}
	\Delta(\theta,r) = \{ z\colon |z| < r, z\ne \psi, |\arg(z - \phi)| > \theta \}
	\qquad\text{for some } r > \phi \text{ and } 0 < \theta < \pi/2.
\end{equation*}
Further assume that when $z\sim \phi$ for  $z\in \Delta(\theta,r)$, $F(z)$ has a singular expansion of the form
\begin{equation*}
	F(z) \sim c\cdot \left( 1 - \frac{z}{\phi} \right)^{-\alpha}
	\qquad\text{for some } c >0 \text{ and } \alpha \notin \{0,-1,-2,\dots\}.
\end{equation*}
Then the coefficients of $F(z)$ satisfy
\begin{equation*}
	[z^n]F(z) \sim \frac{c}{\Gamma(\alpha)} n^{\alpha-1} \phi^{-n}
	\qquad \hbox{as } n \to \infty.
\end{equation*}

We compute first the 242 different roots of $d(z)$ and check that $\rho$ is the only one having modulus $\rho$.
Since $T_b(z)$ is algebraic there exists $r > \rho$ and $0 < \theta < \pi/2$ such that the representation of $T_b(z)$ at $z=0$ admits an analytic continuation to a domain at $z=\rho$ of the form $\Delta(\theta,r)$.
This analytic continuation can be computed from $P(z,T_b)$ using Newton's \textit{polygon algorithm} (see~\cite{FS09}).
This gives an expansion of the \textit{Puiseux} type:
\begin{equation}\label{puisT}
	T_b(z) = t_0 - t_2\left( 1 - \frac{z}{\rho} \right)
	+ t_3\left( 1 - \frac{z}{\rho} \right)^{3/2} + O\left( 1 - \frac{z}{\rho} \right)^2
	\qquad\text{for } z\sim\rho \text{ and  } z\in\Delta(\theta,r),
\end{equation}
where $t_0 = T_b(\rho)\approx 0.000104$ and $t_2\approx 0.002637$.
From there we apply the transfer theorem and obtain the claimed estimate, with $t_3~\approx~0.009747$ and $t = t_3/\Gamma(-3/2) = 3t_3/(4\sqrt{\pi}) \approx 0.00412$.
\qed

\paragraph{Note.}

An alternative for obtaining an upper bound on $\rho$ is to use the class of 3-connected quadrangulations.
This class was first counted by the present authors in~\cite{NRR19,NRR21} as an intermediate step in the enumeration of labelled 4-regular planar graphs.
The radius of convergence of 3-connected quadrangulations counted by number of faces was determined in~\cite{NRR21} as being $\tau = \frac{88 - 12\sqrt{21}}{135}\approx 0.24451$.
Since in a quadrangulation the number of edges is twice the number of faces, the radius of convergence in terms of edges is $\sqrt{\tau}\approx 0.49445$.
Hence $\rho < 1/2$.
This is a better upper bound than $5/8$ but the conclusion regarding the value of $\rho$ remains the same.

%%%%%%%%%%%%%%%%%%%%%%%%%%%%%%%%%%%%%%%%%%%%%%%%%%%%%%%%%%%%%%%%%%%%%%%%%%%%%%%%
%
%	Concluding remarks
%
%%%%%%%%%%%%%%%%%%%%%%%%%%%%%%%%%%%%%%%%%%%%%%%%%%%%%%%%%%%%%%%%%%%%%%%%%%%%%%%%
\section{Concluding remarks}\label{sec:conclude}

Our result could open the way to the enumeration of (labelled) bipartite planar \textit{graphs}, an interesting open problem.
In fact, this was the original motivation for embarking on this project.
For this, one needs the generating functions $T_1(z,\nu,x)$ and $T_2(z,\nu,x)$ of bicoloured 3-connected maps counted additionally according to the number of vertices (marked by the variable $x$).
We have been able to determine the minimal polynomials of both $T_1(z,\nu,x)$ and $T_2(z,\nu,x)$ but they are truly enormous, each containing about $10^7$ monomials in $z$, $\nu$ and $x$.
In~\cite{RW14}, one can find equations, similar to the ones in Section~\ref{sec:enum_bicolor_maps}, relating $T_1(z,\nu,x)$ and $T_2(z,\nu,x)$ to the series of 2-connected and connected bicoloured planar graphs.
Solving these equations (as is done in~\cite{RW14} for series-parallel graphs, a simpler class having no 3-connected graphs, see also~\cite{DRR17} for a discussion comparing with triangle-free graphs), and setting $\nu=0$, one could in principle obtain the series of bipartite planar graphs and deduce a precise asymptotic estimate.
But the task appears daunting given the huge size of the equations and the intricated analysis needed to complete the project in this context.

However, we can solve the simpler problem of counting bipartite \textit{cubic} planar graphs.
Following Tutte, such graphs (or maps) are called \textit{bicubic}.
We use the fact~\cite[Section 11]{T63} that the series of rooted bicubic maps counted by half the number of vertices is precisely equal to $A_b(z)$ given in Equation~\eqref{GF:Mb} (a well-known bijection explains this fact).
If $G(z)$ is the series of 3-connected bicubic maps, then we have (see~\cite{T63}):
\begin{equation*}
	A_b(z) = G\left( z(1 + A_b(z))^3 \right) .
\end{equation*}
As shown in~\cite{MRS74}, the radius of convergence of $G(z)$ is $\tau = 125/512$.
From here one can proceed as in~\cite[Section 3.3]{NRR20}.
The series $D(z)$ of bicubic planar \textit{networks} (essentially edge-rooted 2-connected cubic planar graphs) satisfies the equation
\begin{equation*}
	H(z,D(z)):= D(z) + \frac{z^2}{2}(1 + D(z))- \frac{1}{2}G\left( z^2(1 + D(z))^3 \right) = 0.
\end{equation*}
The radius of convergence $\sigma$ of $D(z)$ is then obtained by solving the system
\begin{equation*}
	\sigma^2(1 + D(\sigma))^3 = \tau, \qquad H(\sigma,D(\sigma))= 0.
\end{equation*}
After checking some rather simple analytic conditions (as in~\cite[Theorem 3]{NRR20}), and observing that a connected bicubic graph is necessarily 2-connected, we obtain the following result:
\begin{theorem}
The number $B_n$ of labelled bicubic planar graphs with $n$ vertices satisfies
\begin{equation*}
	B_n \sim g\, n^{-7/2}\delta^n n!
	\qquad\text{as } n\to\infty,
\end{equation*}
where $g > 0$, $\delta = 1/\sigma\approx 2.035614$, and $\sigma\approx 0.49125$ is the smallest positive root of the irreducible polynomial
\begin{equation*}
	125z^6 + 750z^4 - 4332z^2 + 1000.
\end{equation*}
\end{theorem}

%%%%%%%%%%%%%%%%%%%%%%%%%%%%%%%%%%%%%%%%%%%%%%%%%%%%%%%%%%%%%%%%%%%%%%%%%%%%%%%%
%
%	Acknowledgements
%
%%%%%%%%%%%%%%%%%%%%%%%%%%%%%%%%%%%%%%%%%%%%%%%%%%%%%%%%%%%%%%%%%%%%%%%%%%%%%%%%
\newpage
\section*{Acknowledgements}

We thank the anonymous referee of a first version of this work for having noticed an error in our notion of maps rooted at a bichromatic edge, as well as several suggestions for improving the presentation of the paper.
We are also very grateful to Mireille Bousquet-M\'elou for useful comments, in particular for suggesting the method to obtain the minimal polynomials of $M_1(z,\nu)$ and $M_2(z,\nu)$.

The authors acknowledge the financial support of the Spanish State Research Agency through projects MTM2017-82166-P, PID2020-113082GB-I00, and the Marie Curie RISE research network 'RandNet' MSCA-RISE-2020-101007705.
Additionally, M.N. and J.R. acknowledge support from the Severo Ochoa and Mar\'ia de Maeztu Program for Centers and Units of Excellence (CEX2020-001084-M), and C.R. acknowledges support from the grant Beatriu de Pin\'os BP2019, funded by the H2020 COFUND project No 801370 and AGAUR (the Catalan agency for management of university and research grants).

\bibliographystyle{abbrv}
\bibliography{biblio_bimaps}

\end{document}